\documentclass{conm-p-l}
\usepackage{amssymb}
\usepackage{graphicx}

\usepackage{amsfonts,amsmath, amscd, amsthm, calc, epsfig,epstopdf, ifthen, mathrsfs, mdwlist, pifont,stmaryrd,wasysym,inputenc,}

\usepackage{bbding, comment, float,  microtype, setspace, yfonts}
\usepackage{hyperref,wrapfig,color,}

\newtheorem*{theorem}{Theorem}

\newtheorem*{cor}{Corollary}
\newtheorem*{prop}{Proposition}

\newtheorem*{lemma}{Lemma}

\renewcommand{\(}{\left(}
\renewcommand{\)}{\right)}

\renewcommand{\~}[1]{\overline{#1}}

\renewcommand{\leq}{\leqslant}

\newcommand{\8}{\infty}

\renewcommand{\a}{\alpha}

\newcommand{\Aut}{\text{Aut}\,}

\renewcommand{\Cap}[2]{\underset{#1}{\overset{#2}{\cap} }}

\newcommand{\frakH}{\mathfrak H}

\newcommand{\G}{\Gamma}
\newcommand{\I}{\mathcal{I}}

\renewcommand{\int}{\varint}

\newcommand{\N}{\mathbb{N}}
\renewcommand{\O}{\mathcal{O}}

\newcommand{\Oplus}[2]{\underset{#1}{\overset{#2}{\oplus} }}

\newcommand{\onto}{\twoheadrightarrow}
\renewcommand{\P}{\mathcal{P}}

\newcommand{\Prob}{\textrm{Prob}}

\renewcommand{\qedsymbol}{\includegraphics[width=0.12in]{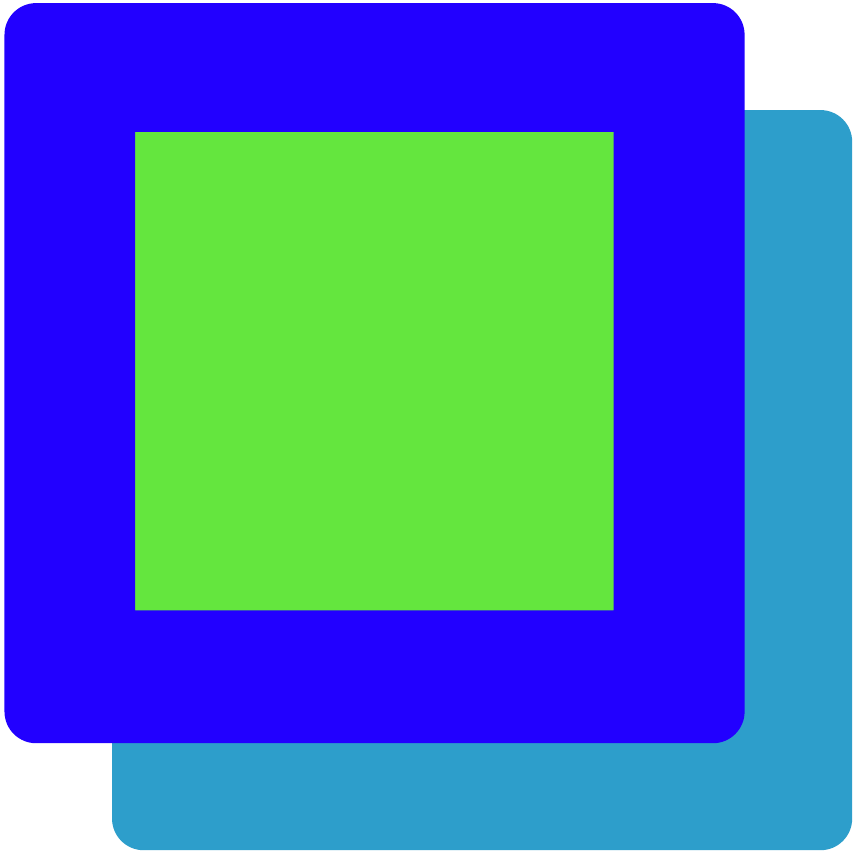}}

\renewcommand{\S}{\mathcal{S}}

\newcommand{\Sqcup}[1]{\underset{#1}{\sqcup}}
\newcommand{\stab}{\mathrm{stab}}
\newcommand{\str}{{\tiny\FiveFlowerOpen}}

\newcommand{\sym}{\mathrm{Sym}}

\newcommand{\Z}{\mathbb Z}

\begin{document}

\title{Le Conte de la Mesure sur les Complexes Cubiques CAT(0)}
\author{Talia Fern\'os}

\address{Department of Mathematics, University of North Carolina, Greensboro}

\email{dr.talia.fernos@gmail.com}

\thanks{The author was partially supported by NSF grant DMS–2005640.}

\subjclass[2020]{Primary 54C40, 14E20; Secondary 46E25, 20C20}
\dedicatory{In honor of Mike Mihalik's 70th birthday, happy birthday!}

\keywords{CAT(0), Cube Complexes, Tits' Alternative}

\subjclass[2020]{Primary 20F65, 20E05, 60B05; Secondary 05C12, 05E18, 06A07,  20F29}

\date{March 10, 2024}

\begin{abstract}
We revisit the topic of probability measures on CAT(0) cube complexes and prove that an amenable group acting on a CAT(0) cube complex, regardless of dimension, necessarily preserves an interval in the Roller compactification. In the finite dimensional case, we prove that there must be an orbit of cardinality $2^N$, where $N$ is bounded by the dimension. This is a slight extension of the author's previuos Tits' Alternative. 
\end{abstract}

\maketitle

\vskip-1in

\begin{figure}
    \includegraphics[width=5in]{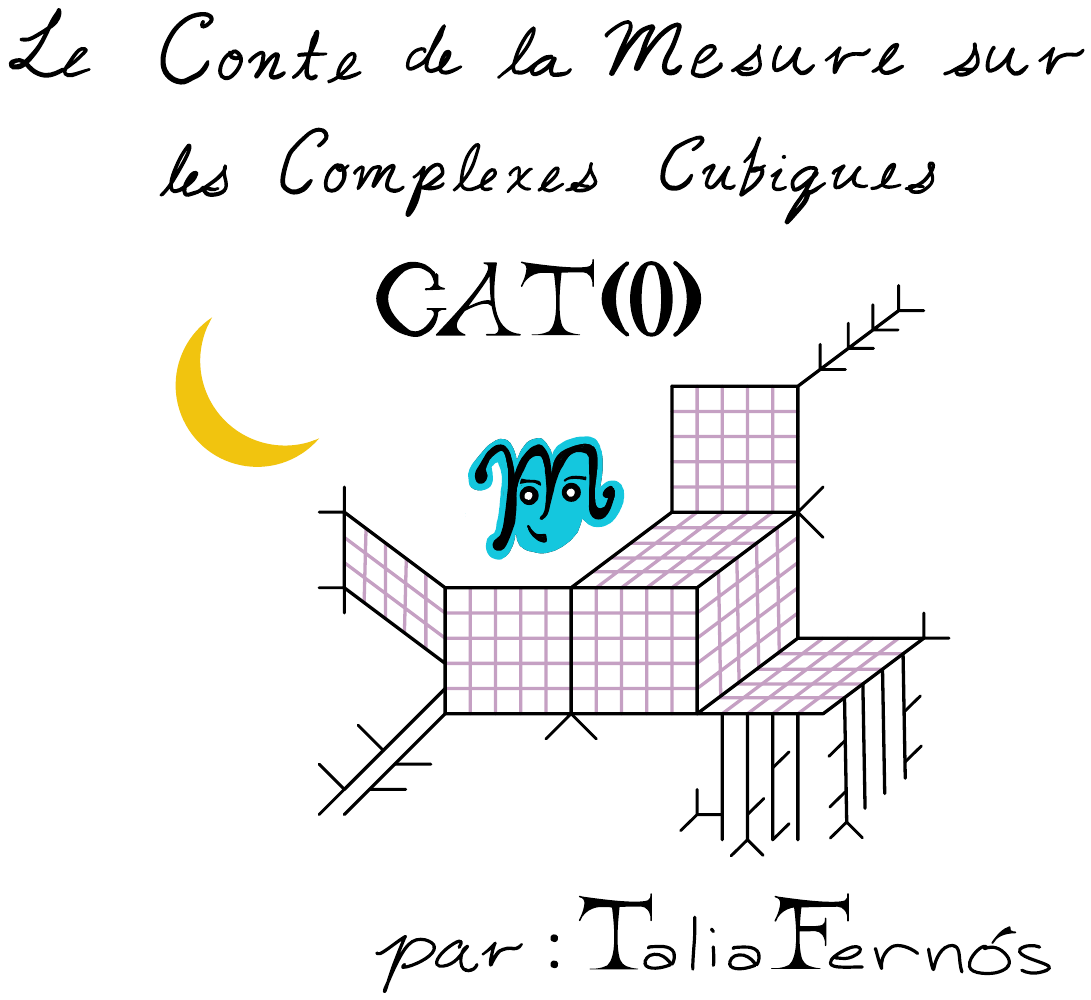}
\end{figure}
\vfill

\specialsection{Introduction}
The idea for this short paper came through preparations for a talk I gave in summer of 2023 at the S\'eminaire Virtuel Francophone -- Groupes et G\'eom\'etrie. The seminar is aimed to be expository in nature and so I thought I would revisit the topic of probability measures on CAT(0) cube complexes, and their associated intervals (see Section \ref{Subsec: Prob Meas}). This association was introduced by Monod and Shalom in \cite[Proposition 4.5]{MonodShalom} for trees and developed further in my joint work with Chatterji and Iozzi \cite{CFI}. It was also used either implicitly or explicitly in \cite{FernosFPCCC, FLM1, FLM2}, and even \cite{FFH}. The key tool for this is the Lifting Decomposition which is discussed in Section \ref{Sec: lifting}.

As I was preparing my lecture, I realized that I could in fact make a small improvement on my previous contribution to the story of Tits' Alternatives. Therefore, I take the opportunity to share that here. I also include some of the drawings I made (using notability on an ipad).

\begin{theorem}
    Let $X$ be a CAT(0) cube complex of finite dimension $D$, and suppose that $\G\to \Aut X$ preserves an interval in the Roller compactification $\~X$. Then $\G$ must have an orbit in $\~X$  of cardinality $2^N$, for some $0\leq N\leq D$.
    
    If $X$ is not finite dimensional and $\G$ is assumed to be amenable then it must preserve an interval $\I\subset \~X$.
    
\end{theorem}

Combining with the Tits' Alternative from \cite{FernosFPCCC}, we obtain the following (the reader should compare statements): 

\begin{cor}[Tits' Alternative]
    Suppose $X$ is a finite dimensional CAT(0) cube complex. Given a group action  $
    \G\to \Aut X$ either $\G$ contains a freely acting free group on 2 generators, or $\G$ has an orbit of cardinality $2^N$ in $\~X$, for some $0\leq N\leq D$, where $D$ is the dimension of $X$. 
\end{cor}

\noindent
{\bf{Examples:}} The following examples show that the bounds above are optimal. We note that all the groups are amenable.

\begin{itemize}
    \item [\str] Consider the standard action of the infinite dihedral group $\mathrm{D}_\8$ on $\Z$. Of course, there is no fixed point but there is an orbit of cardinality 2 in the Roller compactification, namely $\{\pm\8\}$.
    \item [\str] Similarly, $\mathrm{D}_\8^N$ acts on $\Z^N$ and has a unique finite orbit which has  cardinality $2^N$.
    \item [\str]Finally, $\Oplus{n\in \N}{}\mathrm{D}_\8$ acts on $\Oplus{n\in \N}{} \Z$ which is of course infinite dimensional, and does not have a finite orbit in the Roller compactification, but it is an interval. 
\end{itemize}

We note that there are \emph{many} Tits' Alternatives available in the literature \cite{Tits,BestvinaFeighnHandel, SageevWise, CapraceSageev, Fioravanti2018, MartinPrzytycki, OsajdaPrzytycki, GuptaJankiewiczNg, Genevois}.  However, in this brief note, I will not endeavor to give a comprehensive overview of the  Tits' Alternatives, nor to give a thorough account of CAT(0) cube complexes. The reader can see the references provided for further details. 

\noindent
{\bf Acknowledements:} I thank Dafne Sanchez for help with the coloring scheme, inspired by Matisse. I also thank Indira Chatterji, Fran\c cois Dahmani, Anne Lonjou, and Yves Stalder for the invitation to speak at the S\'eminaire Virtuel Francophone -- Groupes et G\'eom\'etrie which inspired me to write this article. I thank all of my collaborators and mentors with whom I have discussed CAT(0) cube complexes for helping me explore how to understand and communicate their beauty. Finally, I thank the National Science Foundation for their generous support of my work, particularly through NSF grant DMS–2005640. 
\specialsection{Basics}

A cube complex is a space obtained by gluing unit cubes isometrically along their faces. Moreover, it is said to be CAT(0) if it is nonpositively curved and simply connected. 

\begin{figure}
    \centering
    \includegraphics[width=0.75\linewidth]{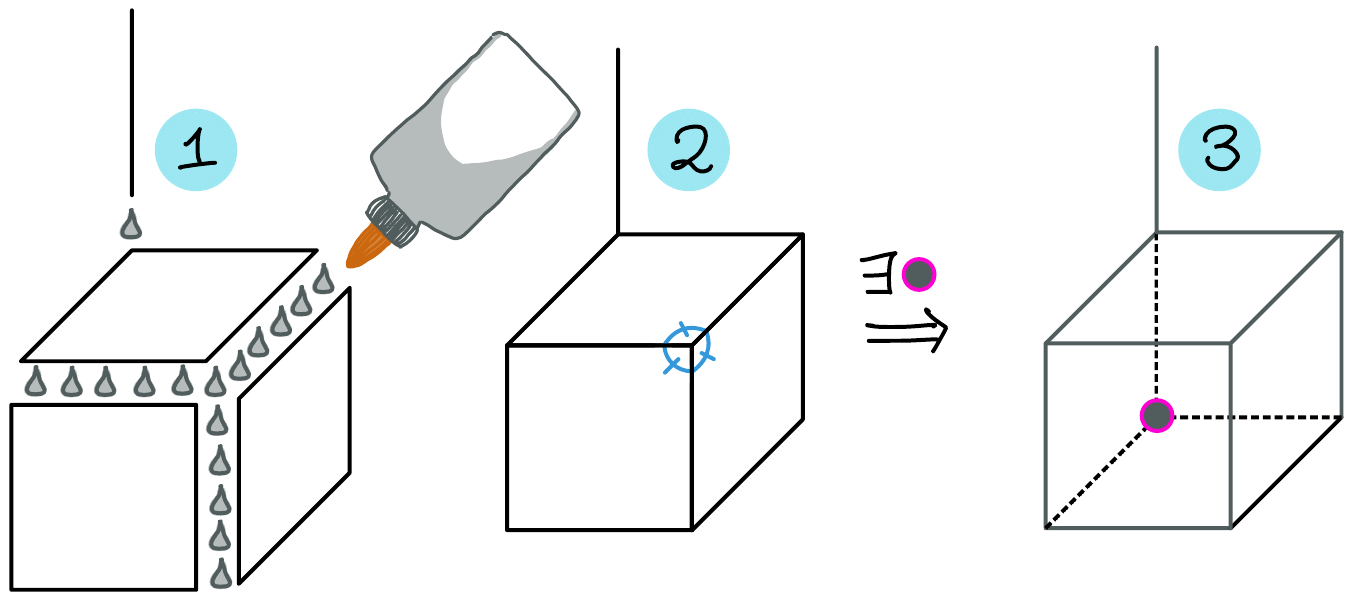}
    \caption{Making a cube complex CAT(0)}
    \label{Fig: Glue}
\end{figure}

You may think of the process of creating a CAT(0) cube complex as gluing together cubes in a particular order. If at some point there is positive curvature that is created locally, it is due to something as in Figure \ref{Fig: Glue}. There, we have created part of the boundary of a  3-dimensional cube by gluing 3 unit squares together. That creates a cone angle of $3\pi/2$. Since $3\pi/2<2\pi$, it is positively curved. However, we may annul the positive curvature by ``filling" in the other side to create a 3-dimensional cube.  Nonpositive curvature of a CAT(0) cube complex is equivalent to the Gromov Link Condition. We refer the reader to \cite{Sageev95, ChatterjiNiblo, Nica, Roller} for more details.

\subsection{Functoriality: Walled Spaces}
A useful way to think of CAT(0) cube complexes is via this functorial construction. Begin with a set $Y$ together with a collection of (nonempty)  two-sided (walls) partitions $\P\subset 2^Y$. Two sided here means that if $h\in \P$ then $Y\setminus h\in \P$. The wall may be thought of as the pair $\{h, Y\setminus h\}$.

\begin{figure}
    \centering
    \includegraphics[width=0.85\linewidth]{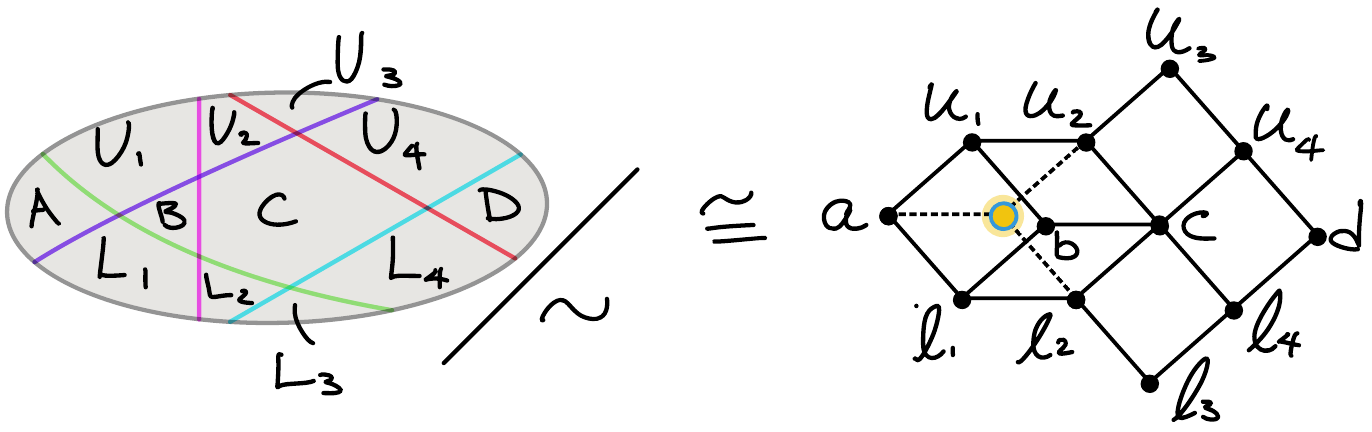}
    \caption{Functorial Construction}
    \label{Fig: Funct Const}
\end{figure}

Next, create a graph by declaring that each maximal non-empty intersection of sets from $\P$ is a vertex. We connect two vertices if their corresponding defining intersections differ by the choice of exactly one side of one partition from $\P$. In Figure \ref{Fig: Funct Const}, you can see a choice of 5 partitions. Each region is not empty. The regions labeled $A$ and $U_1$ differ by the choice of one partition and therefore, their associated vertices $a$ and $u_1$ respectively are connected by an edge.

Once we have done this, we may not have a CAT(0) cube complex. However, we may take the cubical completion (as in Figure \ref{Fig: Funct Const}, where we must add the unlabled ``back" vertex to creat a 3-cube) to annul positive curvature. We note that this can always be done, except that distances may become infinite (if there are infinitely many partitions separating two regions, i.e. the finite interval condition is not satisfied) or the dimension may become infinite (if there are families  of pairwise transverse partitions of unbounded cardinality, see Section \ref{Subs: Comp hs} for the definition of transverse).

Conversely, start with (the vertex set of) a CAT(0) cube complex $X$. Each edge in $X$ belongs to an equivalence class generated by ``being parallel across a square". The compliment of that parallelism class has two sides, and the vertices that belong to each side give the two-sided partition of the vertex set (see Figure \ref{Fig: cut CCC}). 

Finally, if we apply the previous construction to this collection of two-sided partitions, we get  a CAT(0) cube complex, which is cannonically isomorphic to (the vertex set of) $X$.

\begin{figure}\label{Fig: cut CCC}
    \centering
    \includegraphics[width=0.75\linewidth]{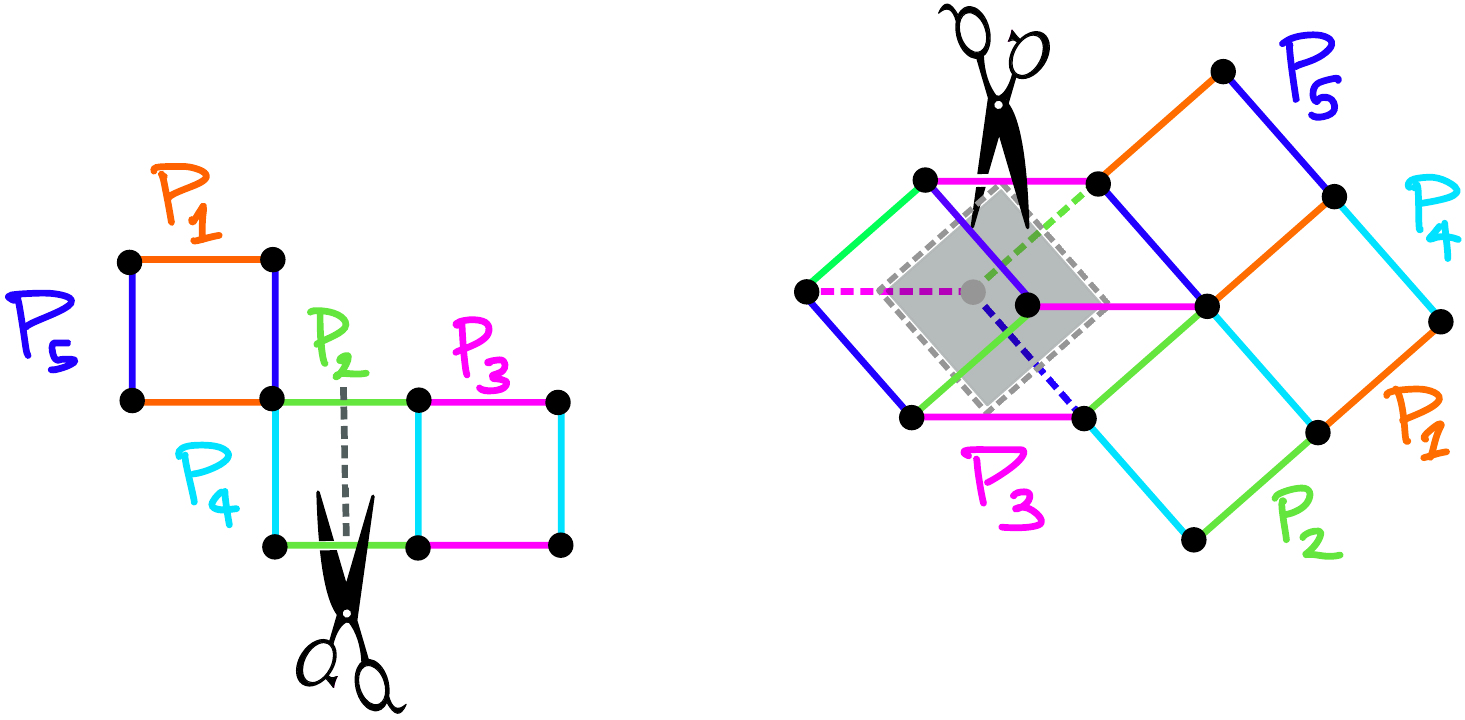}
    \caption{Two CAT(0) Cube Complexes with 10 halfspaces}
    \label{Fig: cut CCC}
\end{figure}

\subsection{The Roller Compactification}\label{Subsec:Roll Comp}
A \emph{halfspace} of a CAT(0) cube complex is one side of a two-sided partition, as discussed above.
We have established that the vertex set, denoted by $X$, is cannonically determined by the halfspace structure   $\frakH\subset 2^X$.  By mapping a vertex in $X$ to the collection of halfspaces that contain it, we obtain an isometric injection:
$$X\hookrightarrow 2^\frakH.$$

Here the extended metric on $2^\frakH$ is given by half the cardinality of the symmetric difference
$d(S,T)= \frac{1}{2}(S\triangle T)\in [0,\8]$. 
While this metric is useful, we will rely on the standard topology on $2^\frakH\cong\mathrm{Map}(\frakH\to \{0,1\})$. This is given by declaring the basic open sets to be cylinder sets, which are themeselves determined by specifying values in \emph{finitely many coordinates}. Equivalently, this is the topology of pointwise convergence of maps $\frakH\to \{0,1\}$. With this topology, $2^\frakH$ is compact. We note that the induced topology may be different than the metric topology on $X$. This is directly comparable to the weak-$*$ topology and respectively the metric topology on a Hilbert space.

Having found an injection of $X$ into the compact space $2^\frakH$, we may take the closure in the image and this defines the \emph{Roller compactification}. By removing the image of $X$ inside this closure, we are left with the \emph{Roller boundary}. 

Note that once we obtain the closure of $X$ in $2^\frakH$, there is a cannonical extension of the halfspaces as partitions of $X$  to $\~X$. Those  are the basic clopen sets for our totally disconnected topology on $\~X$. 

The following is then immediate from the functorial construction, and can be thought of as  forgetfulness. 
\begin{cor}
    Let $X$ be a CAT(0) cube complex with halfspace structure $\frakH$. Let $\frakH'\subset\frakH$ be involution invariant and  $X'$ the associated CAT(0) cube complex. Then the map $2^\frakH\onto2^{\frakH'}$, given by $S\mapsto S\cap \frakH'$ induces a 1-Lipschitz projection $\~X\onto \~{X'}$.
\end{cor}

\begin{figure}
    \centering
    \includegraphics[width=0.5\linewidth]{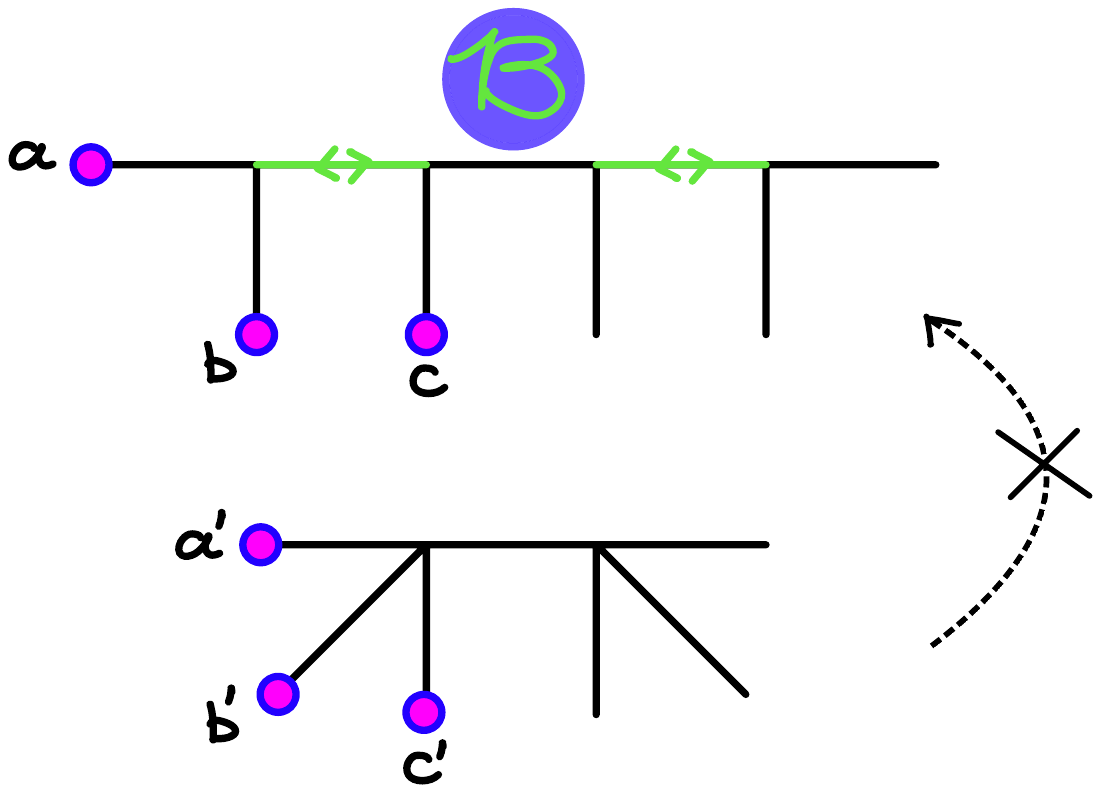}
    \caption{An example where a section to the projection is impossible}
    \label{Fig: Counterex}
\end{figure}

We shall see in Section \ref{Sec: lifting} how and when we can find a cannonical section to this map. The example in Figure \ref{Fig: Counterex} is given by taking $\frakH'= \frakH\setminus \mathcal B$ and a valence consideration shows that an isometric section is impossible. 

\begin{figure}
    \centering
    \includegraphics[width=2in]{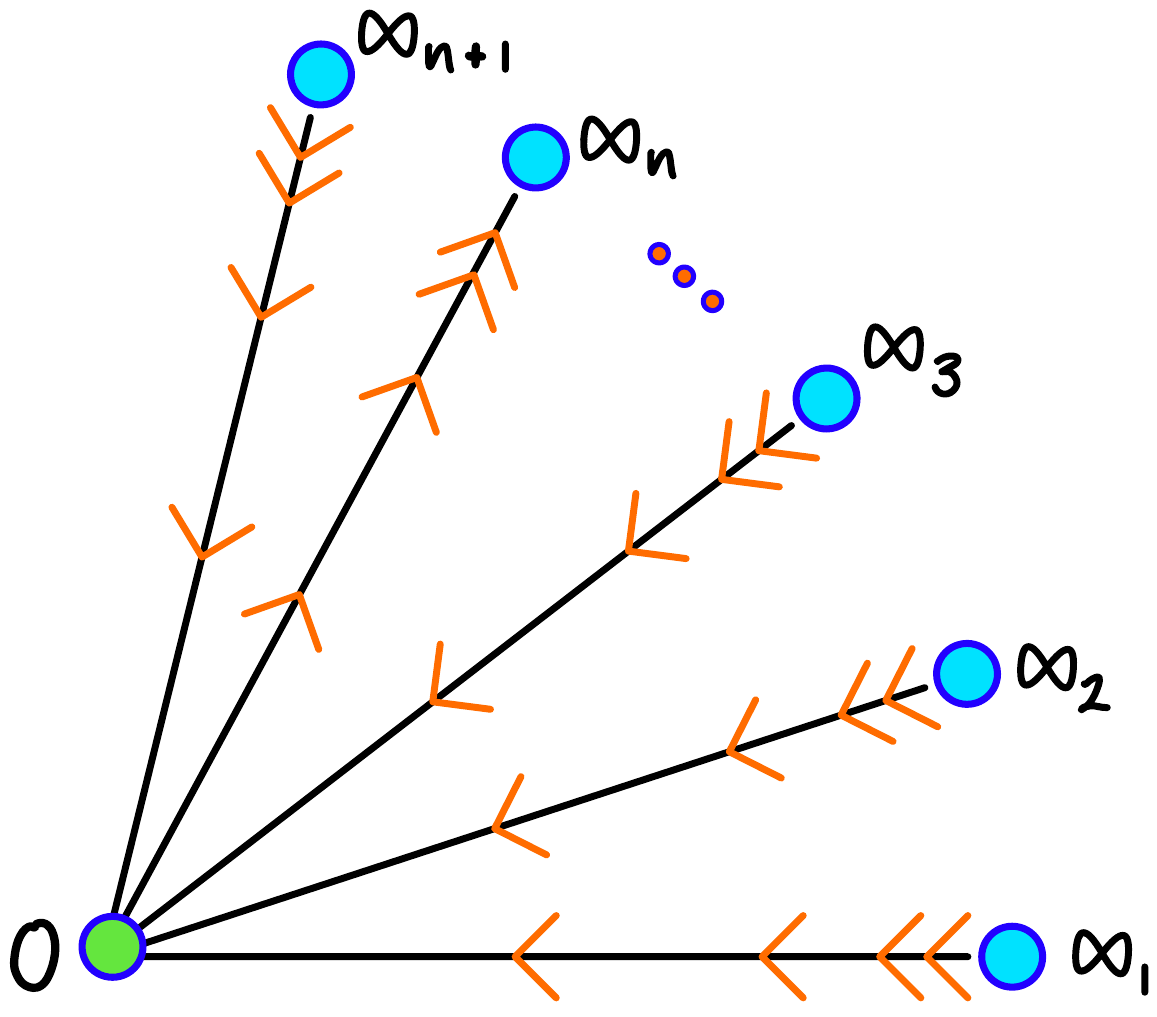}
    \caption{Halfspaces containing $\8_n$.}
    \label{Fig: LocInfinite}
     \vskip.25in
     \end{figure}

The Roller boundary of a tree is, as a set, the visual boundary, but the topology is different when it is not (large-scale)  locally finite. Consider the tree given by identifying at $0$ infinitely many copies of $[0,\8]$, as in Figure \ref{Fig: LocInfinite}. It is not difficult to see that $\8_n\to 0$ since all halfspaces eventually contain $0$.

The Roller compactification  has an important multiplicative property:

$$\~{X_1\times X_2}= \~X_1\times \~X_2.$$
     \begin{figure}
    \centering

     \includegraphics[width=2.5in]{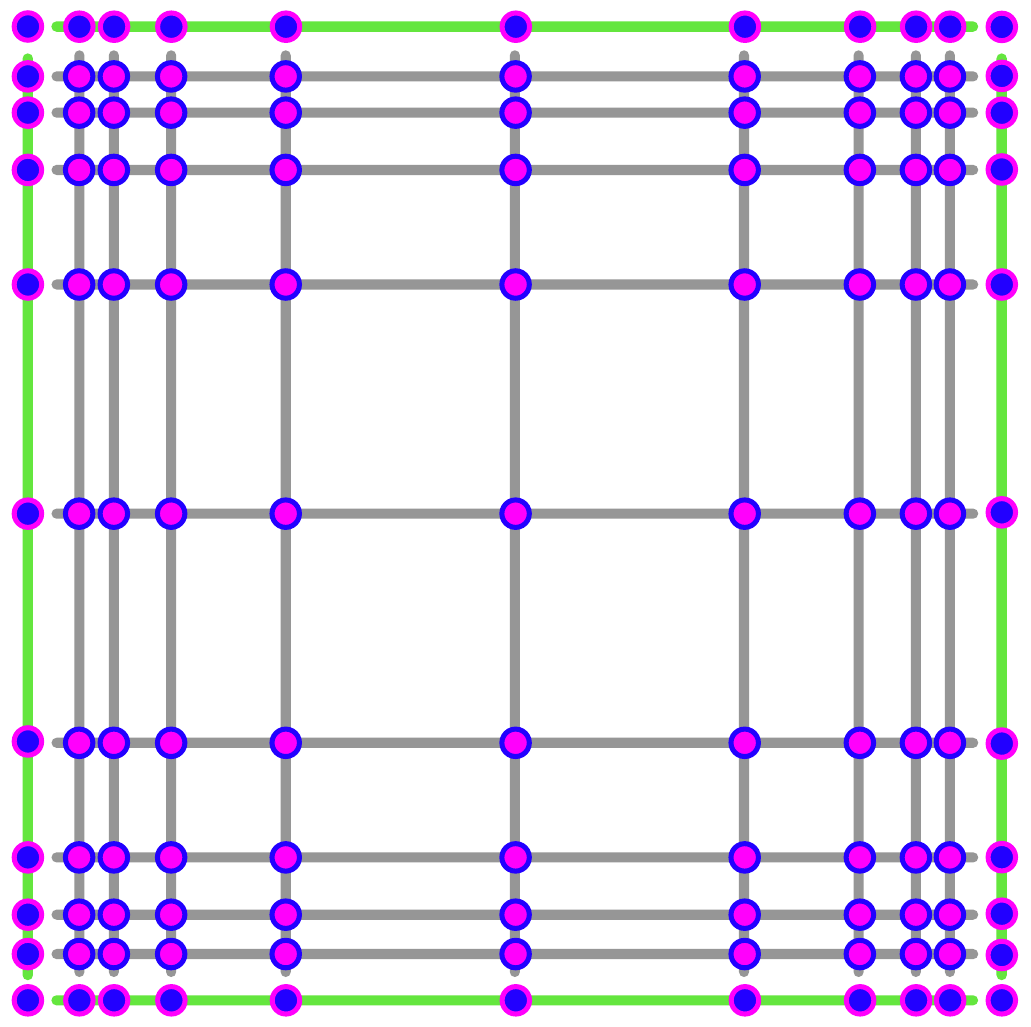}
    \caption{The Roller Compactification for $\Z^2$}
   \label{Fig: Roller Comp Z2}

    \end{figure}
     \begin{figure}
    \centering
    \vskip.25in
    \includegraphics[width=4in]{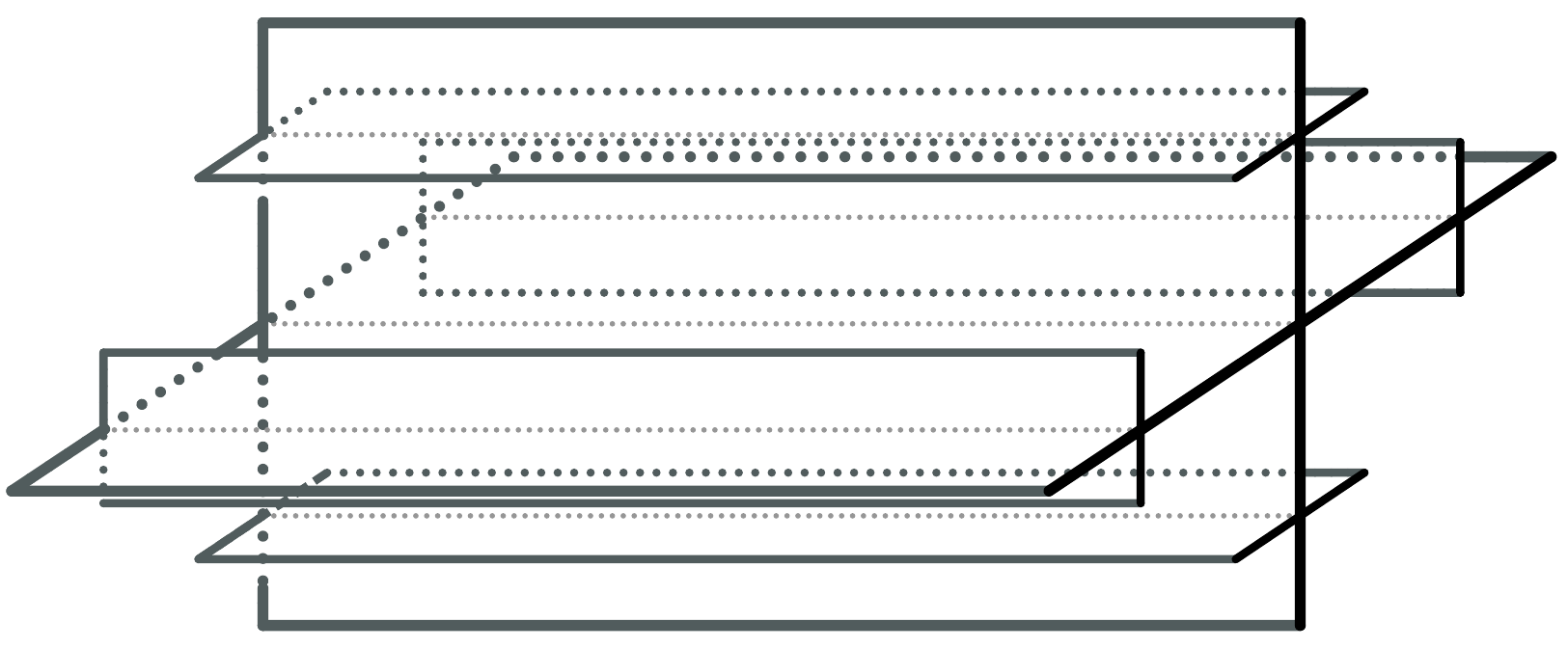}
    \caption{Roller Compactification for $F_2\times \Z$}
    \label{Fig: Roller Comp F2xZ}
\end{figure}

Combining this, with the discussion for trees, we obtain the Roller compactification of $\Z^2$, as in Figure \ref{Fig: Roller Comp Z2}, or for $F_2\times \Z$ as in Figure \ref{Fig: Roller Comp F2xZ}.

\break
\subsection{Comparing halfspaces}\label{Subs: Comp hs}
\begin{figure}
    \centering
    \includegraphics[width=3in]{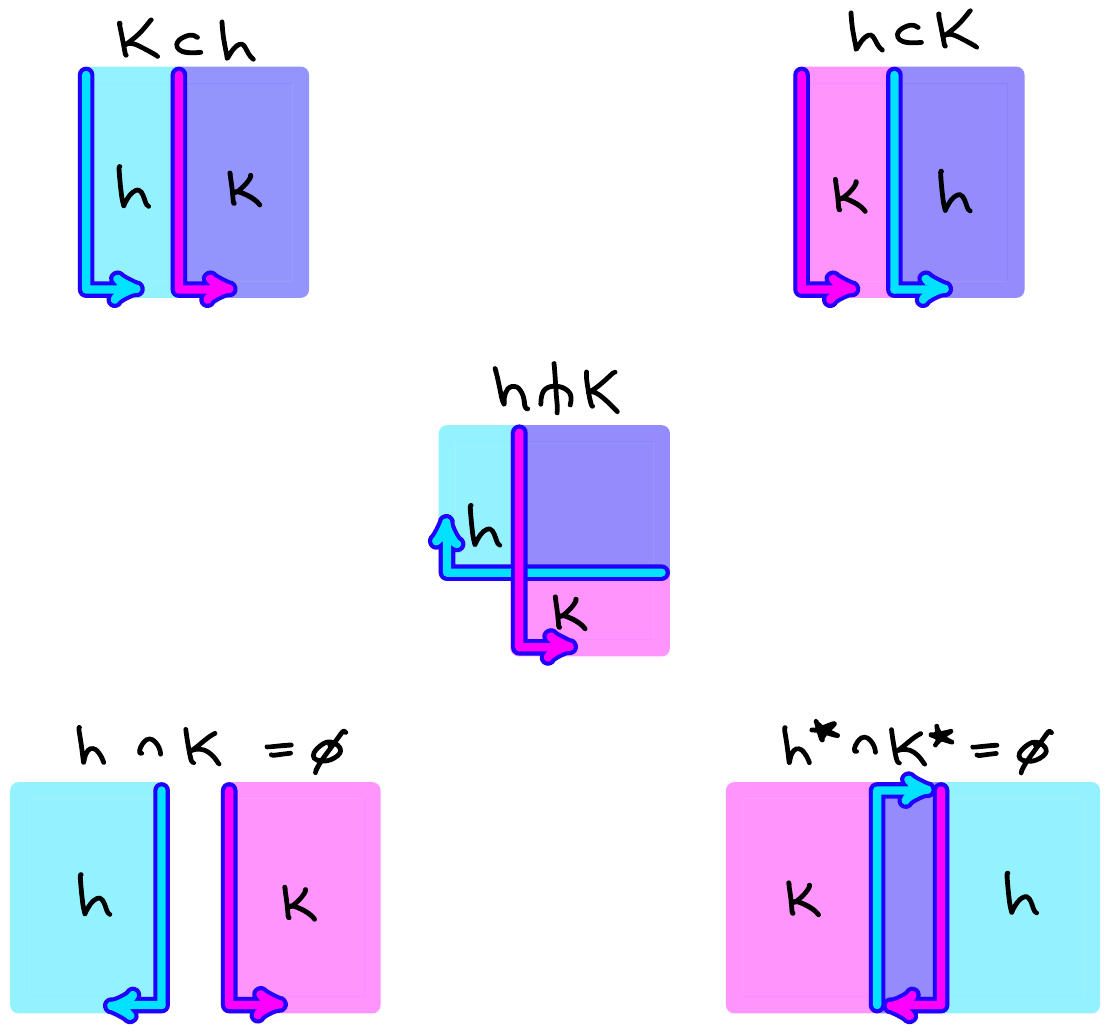}
    \caption{The possible relationships between halfspaces $h,k$, with $k\neq h, h^*$. }
   \label{Fig: Halfspace Pairs}
\end{figure}

Given  a pair of halfspaces $h,k\in \frakH$ such that $k\neq h, h^*$ then either all of the pairwise intersections $h\cap k, h\cap k^*, h^*\cap k, h^*\cap k^*$ are not empty  (and in this case, we say that $h$ and $k$ are \emph{transverse}, and write $h\pitchfork k$) or one of the following other cases hold (see Figure \ref{Fig: Halfspace Pairs}):

\begin{eqnarray*}
    k\subset h \quad\,\,&\qquad\qquad &\quad\,\, h\subset k\\
 h^*\cap k^* =\varnothing &\qquad\qquad & h\cap k =\varnothing.
\end{eqnarray*}

\begin{lemma}
    The dimension of $X$ is the supremum cardinality of subsets $S\subset\frakH$ of pairwise transverse halfspaces.
\end{lemma}

It is straightforward to verify that if $\frakH =\frakH_1\sqcup \frakH_2$ is a nontrivial disjoint decomposition into nonempty involution invariant, pairwise transverse sets and $X$, $X_1$, and $X_2$ are the assocaited CAT(0) cube complexes respectively  then there is a cannonical isomorphism $X\cong X_1\times X_2$. We say that $X$ is irreducible if it is not isomorphic to a product. 

\begin{theorem}\cite[Proposition 2.6]{CapraceSageev} 
    Let $X$ be finite dimensional. Then, there is a cannonical decomposition of $X$ into a product of irreducible CAT(0) cube complexes. 
\end{theorem}

A pair of halfspaces $h, k\in \frakH$, with $h\neq k^*$ are said to be \emph{facing} if $h^*\cap k^*=\varnothing$. This corresponds to the lower right corner in Figure \ref{Fig: Halfspace Pairs}. A triple of halfspaces $h, k, \ell \in \frakH$ are said to be a \emph{facing triple} if they are pairwise facing. 

\break
\subsection{Intervals in $\~X$}\label{Subsect: Intervals}

Let $\I\subset \~X$ be an arbitrary subset. The (involution invariant) collection of halfspaces that separate points in $\I$  is 
$$\frakH_\I:= \{h\in \frakH : h\cap \I, h^*\cap \I\neq \varnothing\}.$$ 

Similarly, we denote the collection of halfspaces that contain $\I$ as $\frakH_\I^+$, and the collection of  compliments of halfspaces that contain $\I$, namely the halfspaces that trivially intersect $\I$ is denoted by $\frakH_\I^-$, meaning that  $\frakH_\I^-=(\frakH_\I^{+})^*$. We have established the following decomposition: 

$$\frakH = \frakH_\I\sqcup \frakH_\I^+ \sqcup \frakH_\I^-.$$

We say that, $\I$ is an \emph{interval} if there exists $x,y\in \~X$ such that  $\I = \Cap{h\in \frakH_{\{x,y\}}^+}{}h$.  In this case we write $\I=\I(x,y)$ and say that it is the \emph{interval between $x$ and $y$}. We also call $x$ and $y$ \emph{endpoints} of $\I$. Clearly, $\frakH_{\I(x,y)}$ does not contain a facing triple.

\begin{figure}
    \centering
    \includegraphics[width=3in]{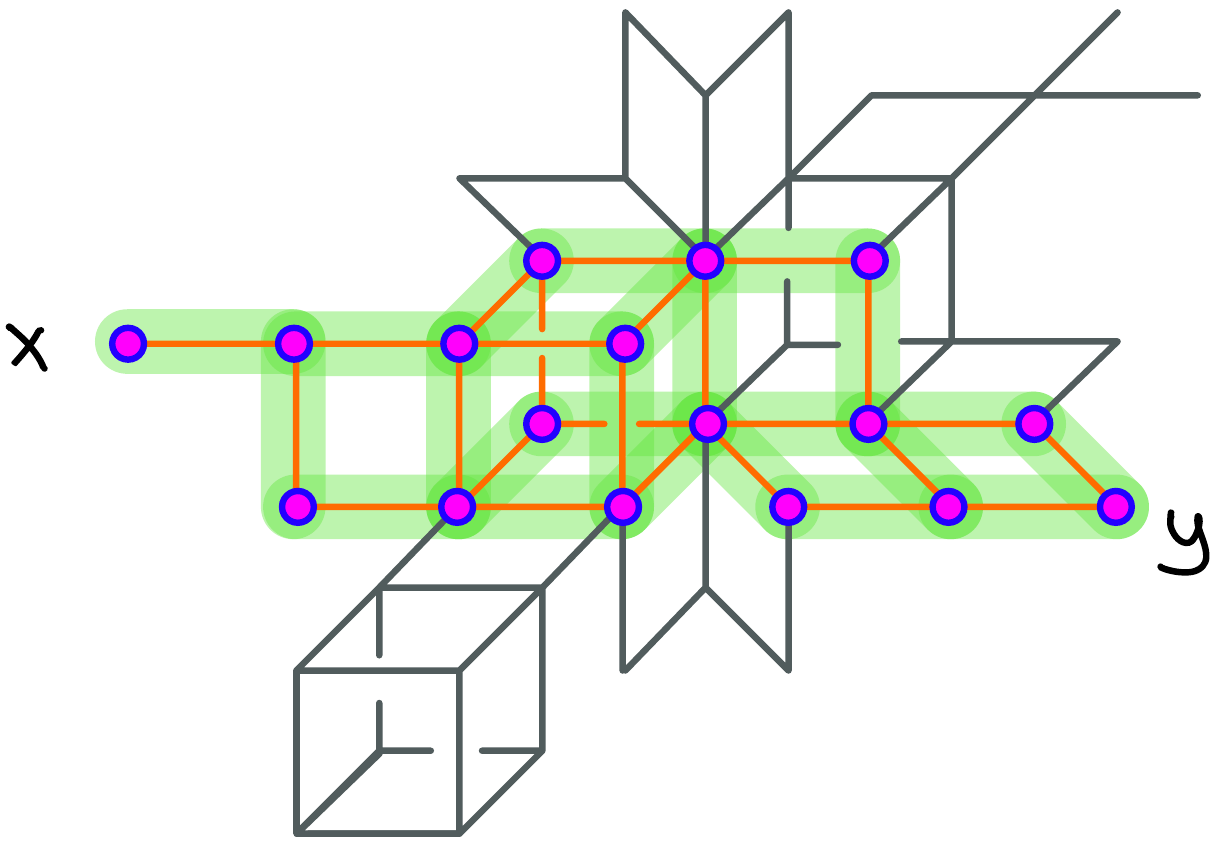}
    \caption{The shaded region is an interval in the ambient CAT(0) cube complex}
    \label{Fig: IntEx}
\end{figure}

\begin{figure}
    \centering
    \includegraphics[width=2.5in]{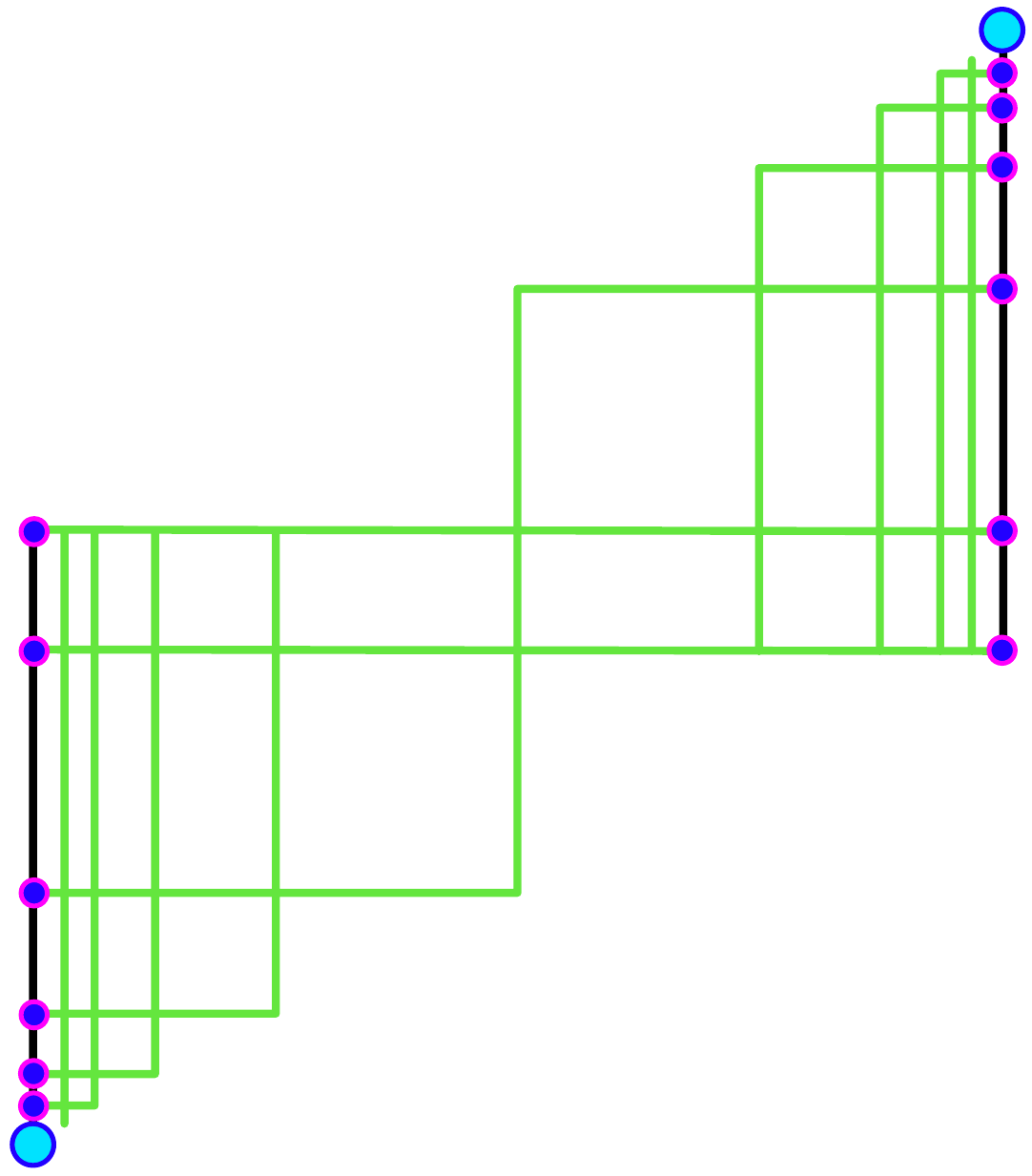}
    \caption{This Roller compactification is an interval}
    \label{Fig: IntEx2}
\end{figure}

The structure of intervals in an arbitrary CAT(0) cube complex can be quite exotic; please see the examples in Figures \ref{Fig: IntEx} and \ref{Fig: IntEx2} or come up with your own!

Nevertheless, the following theorem of \cite{BrodzkiCampbellGuentnerNibloWright}, which relies on Dilworth's Theorem for partially ordered sets, shows that intervals are not too wild. 

\begin{theorem}
    Let $D$ be the dimension of $X$. If $\I\subset \~X$ is an interval then there exists an isometric embedding $\I \hookrightarrow \~{\Z^D}$.
\end{theorem}

\begin{cor}
   If $\I\subset \~X$ is an interval then the set of end points of $\I$ has cardinality $2^N$, for some $0\leq N\leq D$.
\end{cor}

\subsection{The Helly Property in $\~{X}$}\label{Subsec: Helly}

The following is an important property of halfspaces in a CAT(0) cube complex. The first item is the standard Helly property for halfspaces. The second follows from the first by applying the finite intersection property to the collection of  sets $\S$, each of which is clopen in the compact space $\~X$. Comparing it with the classical Helly property for convex sets in euclidean space, we may consider it as telling us that, in a sense, CAT(0) cube complexes are of Helly-dimension 1. 

\break

\begin{lemma}  Suppose that $\S\subset \frakH$ is such that $h\cap k\neq \varnothing$ for all $h,k\in \S$. 
    \begin{itemize}
        \item If $\#\S<\8$ then $\Cap{h\in \S}{}h\cap X\neq \varnothing$.
        \item Otherwise $\~X\supseteq\Cap{h\in \S}{}h\neq \varnothing.$
    \end{itemize}
\end{lemma}

\specialsection{Lifting Decompositions}\label{Sec: lifting}

Let $X$ be a CAT(0) cube complex with associated halfspace structure $\frakH$. We shall say $\S\subset \frakH$ is \emph{consistent} if it satisfies the following two properties:

\begin{enumerate}
    \item If $h\in \S$ then $h^*\notin \S$.
    \item If $h, k \in \frakH$, with $h\subset k$, and $h\in \S$ then $k\in \S$.
\end{enumerate}

We note that, according to the functorial construction, if $\S$ is consistent and $\S\sqcup \S^*= \frakH$ then $\S$ is the collection of halfspaces containing a single point. 

\begin{prop} Let $X$ be a CAT(0) cube complex with associated halfspaces $\frakH$ and $\S\subset \frakH$ consistent. Set $\frakH_\S := 
\frakH\setminus(\S\sqcup \S^*)$ and let $X(\frakH_\S)$ be the associated CAT(0) cube complex. Then the map $2^{\frakH_\S} \to 2^{\frakH}$ given by $E\mapsto E\sqcup \S$ induces an isometric injection 
$\~X(\frakH_\S) \hookrightarrow \~X$, whose image is exactly
$$\Cap{h\in \S}{}\, h\subset \~X.$$

Furthermore, if $\S$ is $\G$-invariant, for some action $\G\to\Aut(X)$ then with the restricted action on 
the image, the above natural injection is  $\G$-equivariant.
\end{prop}

\break
\subsection{Probability Measures}\label{Subsec: Prob Meas}

Consider a probability measure $\mu\in \Prob(\~X)$ and the associated collection of halfspaces $\frakH_\mu^+:= \{h\in \frakH: \mu(h)>1/2\}$. It is straightforward to verify that $\frakH_\mu^+$  is consistent and that the collection of halfspaces $\frakH_\mu=\{h\in \frakH: \mu(h)=1/2\}$ does not contain any facing triples. Applying the Lifting Decomposition to $\frakH_\mu^+$, we get that $\Cap{h\in \frakH^+_\mu}{}\, h$ is isomorphic to $\~X(\frakH_\mu)$.

\begin{lemma}\cite[Lemma 4.7]{CFI} 
    If $\mu\in \Prob(\~X)$ then $\~X(\frakH_\mu)$ is an interval. 
\end{lemma}

\subsection{Medians}\label{Subsect: Medians Intervals}

Let $x,y,z\in \~X$ and consider the associated probability measure $\mu_{(x,y,z)}= \frac{1}{3}(\delta_{x}+\delta_{y}+ \delta_{z})$. A simple parity argument shows  that if $h\in \frakH$ then $\mu_{(x,y,z)}(h)\neq \frac{1}{2}$. Therefore, we have that $\Cap{h\in \frakH_\mu^+}{}h$ is a single point, which we shall call the \emph{median} of the triple and denote it by $m(x,y,z)$.

While this is not the standard definition of the median, it fits nicely within our context. We note that several of the natural properties of the median (e.g. invariance under permutation of the points) are immediate, including the property that $$m(x,x,y)= x.$$

\specialsection{Proof of the Main Theorem}

Suppose $\G$ is amenable acting on $X$. Since $\~X$ is compact and metrizable, there must be a $\G$-invariant probability measure $\mu \in \Prob(\~X)$. By the Lifting Decomposition and the previous lemma, it follows that $\Cap{h\in \frakH_\mu^+}\, h$ is a $\G$-invariant interval.

Suppose now that $X$ is of finite dimension $D$ and that $\G$ is not necessarily amenable but preserves an interval $\I\subset \~X$. By Corollary \ref{Subsect: Intervals} (see also \cite[Corollary 2.9]{FernosFPCCC}), the number of end points on which $\I$ is an interval is $2^{D'}$ for some $0\leq {D'}\leq D$ and we identify these with $\({\Z/2\Z}\)^{D'} \cong \{0,1\}^{D'}$.  By \cite[Proposition 2.6]{CapraceSageev} we have that $\Aut(\{0,1\}^{D'})\cong  \{0,1\}^{D'}\rtimes \sym({D'})$. We may therefore project $\G\to \Aut(\{0,1\}^{D'})$, and without loss of generality, assume $\G\leq \Aut(\{0,1\}^{D'})$.

Let $\G_0 = \G\cap\{0,1\}^{D'}\lhd \G$ and note that $|\G_0|=2^N$ for some $0\leq N\leq {D'}$. Fix  a choice of right $\G_0$-coset representatives $S\subset \G$, with trivial $\{0,1\}^{D'}$-coordinate.  Let $\O\in \{0,1\}^{D'}$ be defined component-wise by $\O_i\equiv 0$, for $i\in \{1, \dots, D'\}$. Note that $\O$ is fixed by $\sym({D'})$ and has trivial stabilizer in $\G_0$. Let $(\O,\a) \in S$. Then $\G_0(\O,\a).\O = \G_0.\O$. We have shown that

$$\G.\O = \Sqcup{(\O,\a)\in S} \G_0(\O,\a).\O = \G_0.\O$$

Therefore, by the orbit stabilizer theorem, we have $$|\G.\O| = |\G_0.\O| = |\G_0|/|\stab_{\G_0}(\O)|= 2^N.$$ 
\hfill \qedsymbol
\break 
\bibliographystyle{amsalpha}
\bibliography{LeConte}

\end{document}